\documentclass{article}
\usepackage[utf8]{inputenc}
\usepackage{amsmath,amsthm,amssymb,amsfonts}

\usepackage{graphicx}
\usepackage{tikz} 
\usepackage{comment}

\newtheorem{thm}{Theorem}

\newtheorem{lem}[thm]{Lemma}

\def\1{{\bf 1}}

\title{Can a rudderless species survive?
}
\begin{document}
\author{
Rinaldo B. Schinazi \footnote{Department of Mathematics, University of Colorado, Colorado Springs, CO 80933-7150, USA;
rinaldo.schinazi@uccs.edu}
}

\maketitle

\begin{abstract}
    Some species of salmon and sea turtle are famously good at finding their birth place to reproduce after having travelled vast expanses of ocean. In contrast, imagine now a species (maybe ancestral to the salmon or turtle) which has to find its birth place to reproduce but has no navigation skills and relies on chance alone. Would such an imaginary species survive? According to our (very simple) model it would survive if and only if the probability that a given individual find its birth place is strictly larger than 1/2. 
\end{abstract}

\section{The model}

 Let $S$ be a countable set. For $x$ and $y$ in $S$, let $p(x,y)$ be the transition probability from $x$ to $y$ for an irreducible discrete time Markov chain {\bf X} on $S$. Let $O$ be a fixed site in $S$. We define a branching Markov chain ${\bf Y}$ as follows. At time $0$, ${\bf Y}$ starts with a single individual at $O$. At every discrete time, if the individual is at $x$ it jumps to $y$ with probability $p(x,y)$ (the transition probabilities of {\bf X}). Before each step the individual has a probability $1-\alpha$ of dying where $\alpha$ is a fixed parameter in $(0,1]$. Whenever the individual returns to $O$ it gives birth to another individual which performs the same dynamics. All individuals behave independently of each other. The process {\bf Y} is said to survive if it has at least one individual somewhere in $S$ at all times. Let $\beta$ be the probability that the Markov chain {\bf X} starting at $O$ eventually returns to $O$. The next result shows that $\beta$ determines whether {\bf Y} may survive.
\begin{thm}
    If $\beta\leq 1/2$ the branching Markov chain {\bf Y} dies out for all $\alpha$ in $(0,1)$.
    If $\beta>1/2$ there exists $\alpha_c\in (0,1)$ such that {\bf Y}  has a positive probability of surviving for $\alpha>\alpha_c$ but dies out for $\alpha\leq \alpha_c$.
\end{thm}

Our branching Markov chain {\bf Y} is a generalization of a model recently introduced by Lebensztayn and Pereira (2023). There, $S=\mathbb{Z}$, $p(x,x+1)=p$ and $p(x,x-1)=1-p$ where $p$ is a parameter in $[0,1]$. In this setting the probability of return is known to be $\beta=1-|1-2p|$, see for instance Grimmett and Stirzaker (2001). Note that $\beta>1/2$ if and only if $1/4<p<3/4$. By direct computation Lebensztayn and Pereira (2023) proved that survival is possible if and only if $p$ is in that range. This note was motivated by the desire to understand their nice result.

As a consequence of our result we see that if the Markov chain {\bf X} is recurrent (i.e. $\beta=1$) then survival is always possible for some $\alpha$. On the other hand if the Markov chain is too transient (i.e. $\beta\leq 1/2$) then survival is possible for no $\alpha$. For instance, survival is possible for the simple symmetric random walk on $S=\mathbb{Z}^d$ for $d=2$ since this is a recurrent chain but not possible for $d\geq 3$, McCrea and Whipple (1940) estimated $\beta$ to be about 0.34 in $d=3$.

Going back to our biological application we can think of $(p(x,y))$ as the probabilities that an individual uses to pick a direction and as $\alpha$ as a measure of the leniency of the environment. Whether the species will survive depends on how likely an individual is to find its birth place in a perfectly lenient environment (i.e. $\alpha=1$). This in turn depends on $S$ and $(p(x,y))$. This model may have some biological relevance in that species with great navigation skills may have evolved from species with poor skills. However, primitive species had to survive in order to evolve into something else.

\section{Proof of Theorem 1}
Following Lebensztayn and Pereira (2023) we define a Bienaym\'e-Galton-Watson process (BGW in short) {\bf Z} that keeps track of the genealogy of the process {\bf Y}. Let $Z_0=1$ and let $Z_1$ be the number of returns of the initial individual to $O$.
Since at each return a new individual is born $Z_1$ also counts the number of children of the initial individual. We can think of $Z_1$ as the number of individuals in the first generation. We define $Z_2$ as the number of children born from the first generation (i.e. the grandchildren of the initial individual) and so on. Since all the individuals are independent of each other and follow the same dynamics ${\bf Z}$ is indeed a BGW process. Moreover, the process ${\bf Z}$ survives if and only if the process {\bf Y} survives.

Note that the total offspring of one individual is the number of times this individual returns to $O$ without being killed. Hence, the mean offspring per individual for the process {\bf Z} is for $0<\alpha<1$,
\begin{equation}
 \mu(\alpha)=\sum_{n\geq 1}\alpha^n p_n(O,O),   
\end{equation}
where $p_n(O,O)$ denotes the probability that the Markov chain {\bf X} starting at time $0$ at $O$ returns to $O$ at time $n$.

We will need the following well known recurrence criterion, see for instance Theorem 1.1 in Chapter 5 in Schinazi (2010). An irreducible Markov chain ${\bf X}$ is recurrent if and only if 
\begin{equation}
    \sum_{n\geq 1} p_n(O,O)=+\infty,
\end{equation}
for some state $O$. We also will need the following result for power series, see Proposition A 1.9 in Port (1994).
\begin{lem}
   Assume that $(b_n)$ is a sequence of positive real numbers such that the series  $\sum_{n\geq 1}b_ns^n$ converges for all $s$ in $[0,1)$. Then,
   $$\lim_{s\to1^-}\sum_{n\geq 1}b_ns^n=\sum_{n\geq 1}b_n,$$
   where both sides of the equality may be infinite.
\end{lem}

There are two cases to consider. Assume first that the Markov chain ${\bf X}$ is recurrent (i.e. $\beta=1$). Then, by Lemma 2 and (2),
$$\lim_{\alpha\to1^-}\mu(\alpha)=\sum_{n\geq 1}p_n(O,O)=+\infty.$$
Since $\mu$ is continuous on $(0,1)$ and $\lim_{\alpha\to 0}\mu(\alpha)=0$,
there exists $\alpha_c$ in $(0,1)$ such that $\mu(\alpha_c)=1$. Since $\mu$ is strictly increasing, $\mu(\alpha)>1$ if and only if  $\alpha>\alpha_c$. Hence, the process ${\bf Z}$ ( and therefore ${\bf Y}$) survives with positive probability if and only if  $\alpha>\alpha_c$. This proves Theorem 1 in the case $\beta=1$.

Consider now the case when the Markov chain ${\bf X}$ is transient. That is, the probability $\beta$ to return to $O$ is strictly less than 1.  By the Markov property, the offspring distribution for the branching process {\bf Z} is for $\alpha=1$,
$$P(Z_1=j|Z_0=1)=(1-\beta)\beta^j,$$
for $j=0,1,2\dots$. Observe that since $0<\beta<1$ this is a proper probability distribution (it is not when $\beta=1$).
Using this offspring distribution we get that the mean offspring $\mu(\alpha)$ for $\alpha=1$ is,
$$\mu(1)=\frac{\beta}{1-\beta}.$$
Note that $\mu(1)>1$ if and only if $\beta>1/2$.
Moreover, $\mu(\alpha)$ can also be expressed using equation (1) for all $\alpha\leq 1$ (including $\alpha=1$).

If $\beta>1/2$ then $\mu(1)>1$. By Lemma 2 the function $\mu$ is continuous on $(0,1]$. It is also strictly increasing. Hence, there exists $\alpha_c<1$ such that $\mu(\alpha_c)=1$ and $\mu(\alpha)>1$ if and only if $\alpha>\alpha_c$. That is, the process ${\bf Y}$ survives with positive probability if and only if $\alpha>\alpha_c$.

On the other hand if $\beta\leq 1/2$ then $\mu(1)\leq 1$. Since $\mu$ is an increasing function, $\mu(\alpha)\leq 1$ for all $\alpha\leq 1$. The process ${\bf Y}$ survives for no value of $\alpha$. This concludes the proof of Theorem 1.


\begin{thebibliography}{100}

\bibitem{Grimmett}

Grimmett, G. R., Stirzaker, D. R. (2001). Probability and Random Processes, 3rd ed. Oxford, NY:
Oxford Univ. Press.

\bibitem{Leben}

Lebensztayn, E., Pereira, V. (2023). On Random Walks with Geometric Lifetimes, The American Mathematical Monthly, DOI: 10.1080/00029890.2023.2274783

\bibitem{McCrea}

McCrea, W. H., Whipple, F. J. W. (1940). Random Paths in Two and Three Dimensions. Proc. Roy. Soc. Edinburgh 60, 281-298.

\bibitem{Port}

Port, S.C. (1994) Theoretical probability for applications. Wiley.

\bibitem{Schinazi}

Schinazi, R.B. (2010). Classical and spatial stochastic processes, second ed. Birkhauser.

\end{thebibliography}
\end{document}